\documentclass[12pt]{article}%
\usepackage[utf8]{inputenc}
\usepackage[all]{xy}

\usepackage{amsmath}
\usepackage{amsfonts}
\usepackage{amssymb}
\usepackage{graphicx}
\providecommand{\U}[1]{\protect\rule{.1in}{.1in}}
\newcommand{\C}{{\mathbb C}}

\newcommand{\R}{{\mathbb R}}

\newcommand{\T}{{\mathbb T}}

\newcommand{\Hh}{{\mathcal H}}

\def\bi{\begin{itemize}}
\def\ei{\end{itemize}}

\usepackage{graphicx}
\usepackage{caption}
\usepackage{subcaption}
\def\com{{[\omega]}}

\newcommand{\ba}{\begin{eqnarray}}
\newcommand{\ea}{\end{eqnarray}}
\newcommand{\bas}{\begin{eqnarray*}}
\newcommand{\eas}{\end{eqnarray*}}
\newcommand{\be}{\begin{equation}}
\newcommand{\ee}{\end{equation}}

\newtheorem{theorem}{Theorem}

\newtheorem{definition}[theorem]{Definition}

\newtheorem{preremark}[theorem]{Remark}

\newtheorem{prenotation}[theorem]{Notation}

\numberwithin{equation}{section}
\numberwithin{theorem}{section}
\begin{document}

\title{{A new approximation method for geodesics on the space of K\"ahler metrics}}
\author{Jos\'e  Mour\~ao\thanks{Department of Mathematics and Center for Mathematical Analysis, Geometry and Dynamical Systems, Instituto Superior T\'ecnico, University of  Lisbon.}, \, Jo\~ao P. Nunes\footnotemark[1] \, and Tom\'as Reis\thanks{Perimeter Institute, Waterloo, Ontario, Canada.}}

\maketitle

%\bigskip

\begin{abstract}

The 
Cauchy problem  for (real analytic) geodesics in the space of K\"ahler
metrics with a fixed cohomology class on a compact  complex
manifold $M$  can be effectively reduced to the 
problem of finding the flow of a related 
Hamiltonian vector field $X_H$, followed by 
analytic continuation of the time to complex time.

This opens the possibility of expressing the geodesic
$\omega_t$ in terms of  Gr\"obner Lie series of the form
$\exp(\sqrt{-1} \, tX_H)(f)$, for local holomorphic functions $f$. The main goal of this paper 
is to use truncated Lie series as a new way of constructing approximate solutions to the geodesic equation.
For the case of an elliptic
curve and $H$ a certain Morse function squared, we approximate
the relevant Lie series by the first twelve terms, calculated with 
the help of Mathematica. This leads to approximate geodesics
which hit the boundary of the space of K\"ahler metrics in finite geodesic time. For quantum mechanical
applications, one is interested also on the non-K\"ahler polarizations
that one obtains by crossing the boundary of the space 
of K\"ahler structures.

\vskip 0.2cm


\end{abstract}

\tableofcontents

\section{Introduction}
K\"ahler manifolds form a rich class of examples where, for example, problems of Riemannian geometry 
lead to very interesting interdisciplinary developments. Notably, the existence of a K\"ahler metric of constant scalar curvature 
on a compact K\"ahler manifold $M$ can be related to algebro-geometric stability properties of $M$. See \cite{Do3} for a recent review. 

If a manifold $M$ has a K\"ahler metric then it has 
an infinite-dimensional space of such metrics. Indeed, let $(M, J_0, \omega)$ 
be a compact K\"ahler manifold with K\"ahler form
$\omega$ and compatible complex structure $J_0$. Then, the space of K\"ahler potentials on $(M, J_0)$
with fixed cohomology class $[\omega]\in H^2(M,\R)$ is naturally identified with an open subset
in the space of smooth real functions on $M$,
$$
\Hh = \left\{\phi \in C^\infty(M) \ :  \ 
\omega_\phi := \omega + i \partial_0 \bar \partial_0 \phi >0\right\} , 
$$ 
where $\partial_0, \overline \partial_0$ denote the $J_0$-Dolbeaut operators. Recall that the positivity condition, $\omega_\phi>0$, means that 
$\omega_\phi$, together with $J_0$, determines a Riemannian metric $\gamma_\phi (\cdot,\cdot) = \omega_\phi(\cdot,J_0\cdot).$

A natural metric on $\Hh$ is the Mabuchi metric 
\cite{M,Sem,Do1,Do2}
$$
g_\phi(h_1, h_2) = \int_M \, h_1 h_2 \, \frac{\omega_\phi^n}{n!},
\qquad h_1, h_2  \in T_\phi \Hh = C^\infty(M).
$$
The expression for the curvature of this metric \cite{Sem, Do1} suggests that the space 
of K\"ahler metrics with cohomology class $[\omega]$,
$$
\Hh_0 = \left\{\omega_\phi \ ,   \phi \in  \Hh\right\} \cong \Hh / \R , 
$$ 
is\footnote{Here, $\R$ acts on ${\mathcal H}$ by adding constants, $\phi \mapsto \phi + a, a\in \R$. The $\partial\bar\partial-$Lemma ensures that $\omega_\phi=\omega_{\phi'}$ if and only if $\phi = \phi' + a$ for some $a\in \R$.}, morally, an infinite-dimensional symmetric space that would correspond to
$$
\Hh_0 \cong G_\C / G,
$$
where $G$ denotes the group of Hamiltonian symplectomorphisms 
of $(M, \omega)$ and $G_\C$ denotes its (non-existent) complexification.

A path of K\"ahler potentials in $\Hh$, $\phi_t$ for $t$ in some open interval in $\R$, is a geodesic  if it satisfies 
$$
\ddot \phi_t = \frac12 \vert\vert \nabla \dot \phi_t\vert\vert_{{\phi_t}}^2.
$$
It is well-known that this is equivalent to the homogeneous complex Monge-Amp\`ere equation 
\begin{equation}\label{ma}
\left(\Omega + i\partial\bar\partial \Phi\right)^{(n+1)}=0,
\end{equation}
where $\Omega$ is the pull-back of $\omega$ to $A\times M$, $\Phi(w,z,\bar z)= \phi_t(z,\bar z)$, 
$\partial$ is the Dolbeaut operator for $(w,z)$, $w$ is an auxiliary complex 
variable on an annulus $A\subset \C$, with $t=\log |w|$ and $(z,\bar z)$ are local coordinates on $M$. As described by Donaldson \cite{Do1}, these geodesic paths 
would correspond to one-parameter subgroups in $G_\C$, generated by ``complexified'' Hamiltonian flows on $M$, with
the Hamiltonian $H$ given by the initial velocity $\dot\phi_0$.  
Hamiltonian evolution in complex time has been studied
both in K\"ahler geometry \cite{Sem, Do1, Do2, BLU, MN}
and in quantum physics \cite{Th, HK, GS, KMN1, KMN2}. Moser's theorem then guarantees that to the geodesic path $\phi_t$ there corresponds a family of diffeomorphisms of $M$,
$\varphi_{it} \in Diff(M)$, which are called Moser diffeomorphisms, relating $\omega$ and $\omega_{\phi_t}$ by\footnote{The fact that these maps are labelled by $it$ instead of 
simply $t$ is explained below.}
$$\varphi_{it}^*(\omega+i\partial_0\bar\partial_0 \phi_t)=\omega.$$ These can be described by a system 
of  non-linear PDEs 
and also in terms of a lifting of the Hamiltonian flow to a complexification of $M$ \cite{BLU}.

In \cite{MN}, in the real-analytic setting and for $M$ compact, it has been shown how the path of Moser diffeomorphisms
can be explicitly described by Hamiltonian evolution of local holomorphic coordinates, analitically 
continued to complex time. For sufficiently small $|t|$, 
one defines a new global complex structure $J_{it}$ on $M$ 
(which is biholomorphic to $J_0$) via new local ($J_{it}$-holomorphic) coordinates defined by the Lie series
$$
z_{it} = e^{it X_H} z,
$$
where $X_H$ is the Hamiltonian vector field of $H$ with respect to $\omega$. 
One then obtains a geodesic path of K\"ahler structures $(M,\omega_{\phi_t},J_0)$. In the symplectic description, 
which we will use below, one fixes the symplectic form rather than the complex structure 
(see Theorem 4.1 and Proposition 9.1 in \cite{MN}), so that the geodesic path becomes
$$
(\omega, J_{it})= \varphi_{it}^*(\omega_{\phi_t},J_0). 
$$

The standard approach to the problem of existence of geodesics for the Mabuchi metric, 
with a given regularity, is based on the so-called continuity method for the complex homogeneous 
Monge-Amp\`ere equation. (See \cite{Do3} and references therein.) In particular, Rubinstein and Zelditch show in \cite{RZ1,RZ2,RZ3} 
that the Cauchy problem for geodesics in the Mabuchi metric, with prescribed initial point and velocity, is ill-posed in general for $C^3$ metrics.  
In the real-analytic context that we consider in this work, the Cauchy problem admits solutions for a short time \cite{M,Sem,Do3,RZ1,RZ2,RZ3}.

In this paper, following the ideas in \cite{MN}, we wish to propose a different method to
study this problem by
reducing the Cauchy problem for geodesics in $\Hh_\com$ to finding the associated 
$\omega$-Hamiltonian flow followed by an appropriate complexification, in the setting of 
the Gr\"obner theory of Lie series of vector fields
\cite{Gro}.
The main goal of the present paper consists, then, in an initial numerical exploration, in a relatively simple 
concrete geometric situation, of this method of construction of approximate 
solutions to the geodesic equations for the Mabuchi metric. The approximation scheme will consist in  
taking only the first $N$ terms in the relevant 
Lie series. We hope that this will stimulate future work, both analytical and numerical, on this approach to 
the study of geodesics for the Mabuchi metric, different from the one based on the Monge-Amp\`ere equation. 
We note that this corresponds effectively, for analytic geodesics, to replacing the non-linear PDE (\ref{ma}) by the system of ODEs corresponding to 
the flow of $X_H$.

We consider the Cauchy problem with initial flat K\"ahler metric
 on the two--dimensional
torus $\T^2$ and $\frac {\partial \phi_t}{\partial t}|_{t=0}=H$, with $H$
the square of a particular Morse function on $\T^2$.  With the help of Mathematica and of 
the supercomputer Baltasar from CENTRA/IST, we calculate approximately ($N=12$), for different values 
of $t$, the conformal factor of the metrics along the geodesics. The solution 
remains K\"ahler for geodesic time $t$ inside the interval 
$t \in (-T_1, T_2)$ for certain finite positive values $T_1, T_2$. It hits 
the boundary of the space of K\"ahler metrics both 
at (negative) time $-T_1$ and (positive) time $T_2$.
Indeed, for $t: \, t > T_2$ and $t < -T_1$,
the solution describes open regions with positive metric as well as open regions with negative definite ``metric''.
While from the strict point of view of Riemannian geometry these solutions are anomalous, they are still very interesting from 
the point of view of the (geometric) quantization of the underlying manifold, as they 
correspond to so-called mixed polarizations (see \cite{KMN2} for a discussion in the context of toric manifolds). Therefore, 
we exhibit the behaviour of the metric also in this region.

We remark that Lie series have been also 
successfully applied to the approximate integration of 
ordinary differential equations, in particular in celestial mechanics
(see  eg \cite{BHT, ED, HLW}).

\section{Geodesic equation on the space of K\"ahler metrics and imaginary time Hamiltonian symplectomorphisms}
\label{cst}

For the compact K\"ahler manifold $(M, J_0, \omega)$ and in the context described above, consider the following 
Cauchy problem for geodesics for the Mabuchi metric,
\be
\label{e-cp}
\left\{
\begin{array}{rcl}
\ddot {\phi_t} &=& \frac 12 ||\nabla \dot \phi_t||_{\phi_t}^2 . \\
\phi_0 &=& 0 \\
\dot \phi_t |_{t=0} &=& H \,  
\end{array}
\right.
\ee
with real analytic initial data.

In \cite{MN}, the following algorithm has been proposed for reducing the 
highly nonlinear partial differential equation 
(\ref{e-cp}) (which is equivalent to the HCMA equation \cite{Sem, Do1, Do2}) to finding the flow of  the Hamiltonian vector field 
$X_H$ of $H$, for the initial symplectic form $\omega$, followed by
analytic continuation to imaginary time.

\begin{itemize}
\item{\bf Step 1: Lie series of $J_0$-holomorphic coordinate functions} -- Let $U_\alpha$ be an open cover of $M$ and $(U_\alpha, z^\alpha_1, \cdots , z^\alpha_n)$, denote 
$J_0$-holomorphic coordinate charts.
For $\tau \in \C$, and for sufficiently small $|\tau|<T$, for some $T>0$, find the Lie series of $X_H$ 
for every $J_0$--holomorphic coordinate function and 
\be
\label{e-gls}
 z^\alpha_j (\tau) = e^{\tau X_H}(z^\alpha_j) = \sum_{k=1}^\infty \frac{\tau^k}{k!} \, (X_H)^k \, (z^\alpha_j). 
\ee

\item{\bf Step 2: Define the Moser isotopy} -- In this step, one uses the 
constructive proof of Theorem 2.6 of \cite{MN}  to turn
the complex symplectomorphisms (\ref{e-gls}) into Moser 
diffeomorphisms $\varphi_\tau$ such that
\be
\label{e-gls2}
z^\alpha_j (\tau) = \varphi_\tau^*(z_j^\alpha) =  e^{\tau X_H}(z^\alpha_j), \quad |\tau | < T . 
\ee
The functions $z^\alpha_j (\tau)$ are $J_\tau$--holomorphic for the 
complex structure 
$$J_\tau = \varphi_\tau^* (J_0) :=
(\varphi_{\tau *})^{-1} \circ J_0 \circ \varphi_{\tau *} $$

\item{\bf Step 3: Restricting to imaginary time and geodesics} --
As shown in Proposition 9.1 of \cite{MN},
by restricting the Moser isotopy of step 2, $\varphi_\tau$,  to imaginary
$\tau = i t, t \in \R$, the path of symplectic forms
\be
\label{e-kiso}
\omega_t = (\varphi^{-1}_{it})^* (\omega),   \quad |t| < T,
\ee
is a geodesic path in $\Hh_\com$ and its K\"ahler
potential is a solution of the Cauchy problem (\ref{e-cp}).
The expression for the K\"ahler potential 
in terms of the imaginary time syplectomorphisms 
is given by (4.1)--(4.3) and (6.7) of \cite{MN}.
Below, however, we will focus on finding approximate  
expressions for the geodesics in terms of the K\"ahler forms
as in (\ref{e-kiso}).

\end{itemize}

\section{New approximate method for finding geodesics and description of the computational method}
\label{ss-ca}

% In order to compute the integral curves for the Hamiltonian flow, one must only compute a finite number of derivative of the components of the Hamiltonian vector field, which can, generally, be easily obtained analytically. Were we interested only on the integral curves themselves, we might have computed this derivatives through some finite-differences method. However, the purpose of computing the \textit{real time flows} is continue them analytically in $t$ in order to obtain the complex time flow, 
% and that cannot be done with a purely numerical approximation of the flows. 

% However, the expression (\ref{e-gls}) allows just as well for a symbolical approximation of the integral curves, which can then be analytically extended. Furthermore, this expression can then be further differentiated symbolically, thus avoiding the enormous source of numerical error that is the finite differences method. Thus one can follow the analytical procedure outlined in the previous section to obtain the conformal factor from the integral curves on the symbolical expression and thus obtain a symbolical approximate expression for the factor itself without adding any further source of numerical error between the truncation of the Lie series and the final expression. (This is not completely true in practice, as computing the conformal factor requires a quadratic operation on the spatial derivatives of the flow, after which one has to discard the resulting $\mathcal{O}(t^{N+1})$ terms).

%%%%%%%%%%%%%%% texto Ze  %%%%%%%%%%%%%%%%%%%%%%%%%%%%%%%%%%%

Let the geodesic $\omega_t$ in (\ref{e-kiso}) be written in the form 
$$
\omega_t = \frac i2 \sum_{j, \bar j} h_{j \bar j}(z, \bar z; t) \, dz_j \wedge d\bar z_{\bar j} .
$$
Then (\ref{e-kiso}) is equivalent to,
\be
\label{e-kiso2}
(\varphi_{it})^* (\omega_t) =  \omega,   \quad |t| < T,
\ee
or, in local coordinates,
\be
\nonumber
 \sum_{j, \bar j} h_{j \bar j}(z(it), \bar z(it); t) \, dz_j(it) \wedge d\bar z_{\bar j}(it)
=  \sum_{j, \bar j} h_{j \bar j}(z, \bar z; 0) \, dz_j \wedge d\bar z_{\bar j}.
\ee
For complex one--dimensional manifolds one obtains
\be
\label{e-cft}
h_{1 \bar 1}(z, \bar z,t) =  \frac {{h_{1 \bar 1}(z, \bar z, 0)}}{\frac{\partial z(it)}{\partial z} \frac {\partial \bar z(it)}{\partial \bar z}-
\frac{\partial z(it)}{\partial \bar z} \frac {\partial \bar z(it)}{\partial z}}.
\ee

We define the $N$-th order approximation to the Lie series in
(\ref{e-gls}) as
\be
\label{e-gls3}
 z_j (\tau; N) =  \sum_{k=0}^N \frac{\tau^k}{k!} \, (X_H)^k \, (z_j)   \, .
\ee
\begin{definition}
The Lie series $N$--th order approximation to the geodesics is 
 (for $n=1$) defined to be the path of metrics given by the 
following conformal factors
\be
\label{e-cft2}
h_{1 \bar 1}^{(N)}(t) =  \frac {h_{1 \bar 1}(z, \bar z, 0)}{\frac{\partial z(it; N)}{\partial z} \frac {\partial \bar z(it; N)}{\partial \bar z}-
\frac{\partial z(it; N)}{\partial \bar z} \frac {\partial \bar z(it; N)}{\partial z}} = \sum_{k=0}^N a_k(x,y)t^k. 
\ee
\end{definition}

The computational implementation is then straightforward. We used Mathematica, as it provides the necessary tools for heavy 
symbolic manipulations. Nevertheless, some implementation challenges became quickly evident. 
Given a reasonably non-trivial Hamiltonian, the successive application of the chain and Leibniz's rules leads to an exponential growth of the number of additive terms in each of the terms of the series. This means that for a reasonable approximation of $N>4$ the final formula is quite long. 
We thus resorted to the \textit{Baltasar} supercomputer, run by CENTRA, and used a node of 48 cores to run the code.

%\vskip 4cm
\section{Description of the results.}
\label{s-3}

\begin{figure}
    \centering
    \begin{subfigure}[b]{0.4\textwidth}
        \includegraphics[width=\textwidth]{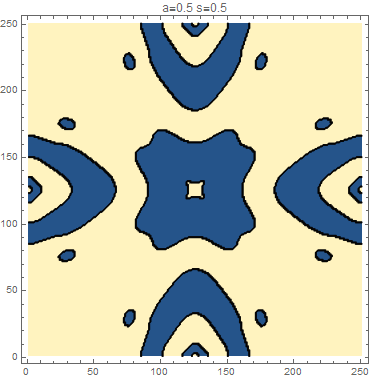}
        \caption{$t=0.5$}
        \label{figposhd}
    \end{subfigure}
    ~ %add desired spacing between images, e. g. ~, \quad, \qquad, \hfill etc. 
    \begin{subfigure}[b]{0.4\textwidth}
        \includegraphics[width=\textwidth]{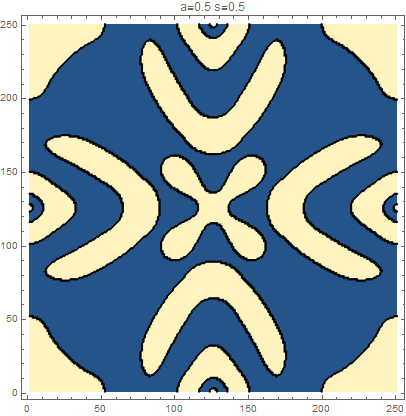}
        \caption{$t=-0.5$}
        \label{figneghd}
    \end{subfigure}
    \caption{Evolution of the sign of $h_{1 \bar 1}$. Regions with positive sign $(+)$ of the conformal factor are in tan and the blue is on regions with 
negative sign $(-)$.}
%\label{fig:animals}
    \label{fighd}
\end{figure}

For the numerical analysis we took the Hamiltonian:
\begin{equation}
H = \frac{1}{8} \left(\sin^2(\pi x) + \sin^2 (\pi y)\right)^2
\end{equation}
in the unit square $[0,1]\times[0,1]$ with initial K\"ahler structure given by the flat metric and the standard complex structure with local holomorphic coordinate, $z=x+iy$. 
An approximate expression for the conformal factor was calculated using 
12 terms of the Lie series.  Except for Figure \ref{fighd}, all the graphs where obtained by sampling 
points on a uniform $50\times 50$ lattice (for each of the the subfigures in Figure \ref{fighd} a $200\times 200$ 
lattice was used instead).

\begin{figure}[!ht]
  \centering
    \includegraphics[width=0.6\textwidth]{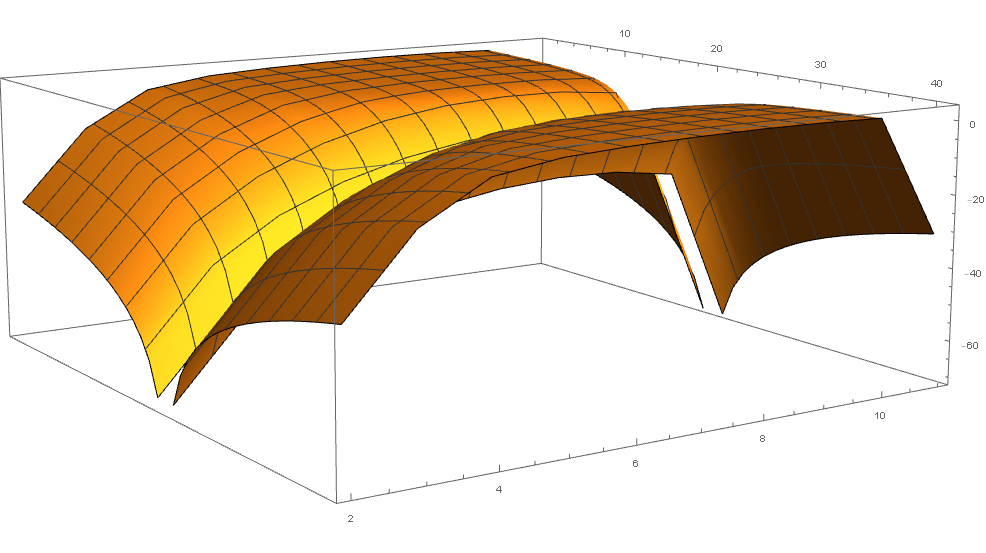}
      \caption{Plot of $log\left(\left|a_{12}t^{12}/h_{1\bar1}^{(11)}(t)\right|\right)$ along the diagonal $x=y \in [0, 0.5]$, 
where $a_j$ is $j$th term in the series and $t$ a value in (imaginary) time. (See equation (\ref{e-cft2}).) 
The valley in the figure corresponds to $t=0$ and the time axis runs in the transverse direction 
for $t\in [-1,1]$.}
%The axis on the bottom facing the reader holds the value of $l\in [0,0.5]$ while the remaining axis cover $s\in [-1,1%]$
\label{figerr}
\end{figure}

%We begin with a brief comment on error analysis and on the chosen intervals. 
It is known  that below a certain 
value of geodesic time the Lie series is absolutely convergent \cite{Gro}. In our analysis, to estimate the error
of the truncation of the series,  we used the ratio of the absolute value 
of the last term considered in (\ref{e-cft2}) to the absolute value of the 
sum of all the lower order terms (see the related discussion around equation
(66) of \cite{ED}). This gives us some indication 
of convergence and of the magnitude of the error.

\begin{figure}
    \centering
    \begin{subfigure}[b]{0.3\textwidth}
        \includegraphics[width=\textwidth]{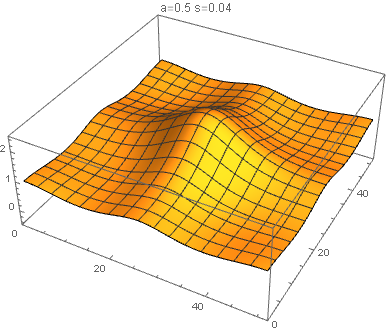}
        \caption{$t=0.04$}
        \label{fig:lp1}
    \end{subfigure}
    ~ %add desired spacing between images, e. g. ~, \quad, \qquad, \hfill etc. 
    \begin{subfigure}[b]{0.3\textwidth}
        \includegraphics[width=\textwidth]{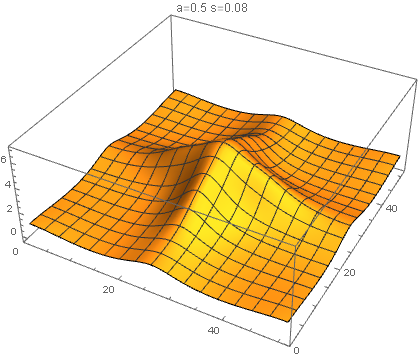}
        \caption{$t=0.08$}
        \label{fig:lp2}
    \end{subfigure}
    ~ %add desired spacing between images, e. g. ~, \quad, \qquad, \hfill etc. 
    \begin{subfigure}[b]{0.3\textwidth}
        \includegraphics[width=\textwidth]{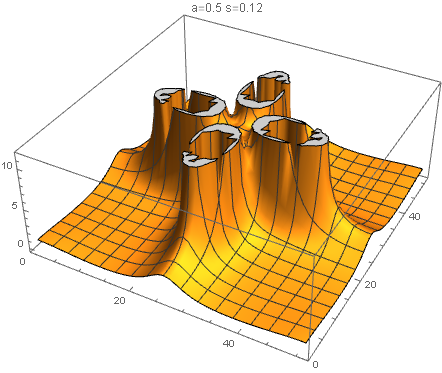}
        \caption{$t=0.12$}
        \label{fig:lp3}
    \end{subfigure}
    \caption{Positive small imaginary time evolution of $h_{1 \bar 1}$}
%\label{fig:animals}
    \label{figlowpos}
\end{figure}
In Figure \ref{figerr}, we consider the error 
%plot the logarithm of the absolute value of the ratio bewteen the last term considered in the sum in (\ref{e-cft2}) and the 
%sum of the lower order terms
for points in the first half of the diagonal 
of the unit square (this is $(x,y)$ such that $x=y<0.5$) and for 
values of imaginary time $t$ in $[-1,1]$. 
We observe that for $t<0.5$ we have negative values of the error indicator for all sampled points, while 
for larger $t$ the same is not true and inclusively we find positive values of the logarithm that suggest a 
significant error. As such we restricted our analysis to $t$ smaller than $0.5$.

Let us begin the analysis of the conformal factor with $t>0$ in Figure \ref{figlowpos}. In Figure \ref{fig:lp3}, 
there are 
regions where the metric is no longer positive definite. (For example, in the regions around the saddle points and the maximum point of $H$. 
See also Figure \ref{fighd}.)  
Although the metric is no longer everywhere 
positive beyond this value of time, it is still interesting for applications in geometric quantization. The critical time, that is 
the earliest time at which the conformal factor is not strictly positive, is higher than but close to $t= 0.118$.

\begin{figure}
    \centering
    \begin{subfigure}[b]{0.3\textwidth}
        \includegraphics[width=\textwidth]{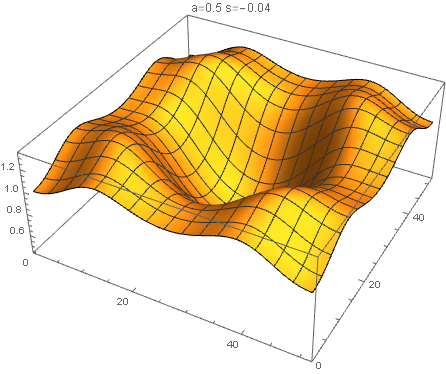}
        \caption{$t=-0.04$}
        \label{fig:ln1}
    \end{subfigure}
    ~ %add desired spacing between images, e. g. ~, \quad, \qquad, \hfill etc. 
    \begin{subfigure}[b]{0.3\textwidth}
        \includegraphics[width=\textwidth]{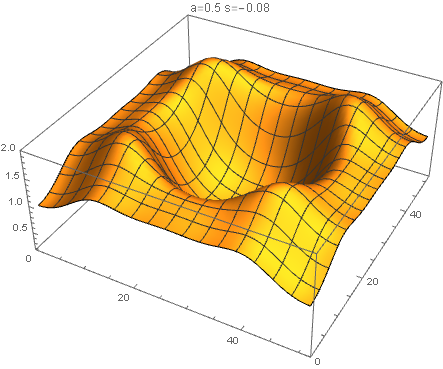}
        \caption{$t=-0.08$}
        \label{fig:ln2}
    \end{subfigure}
    ~ %add desired spacing between images, e. g. ~, \quad, \qquad, \hfill etc. 
    \begin{subfigure}[b]{0.3\textwidth}
        \includegraphics[width=\textwidth]{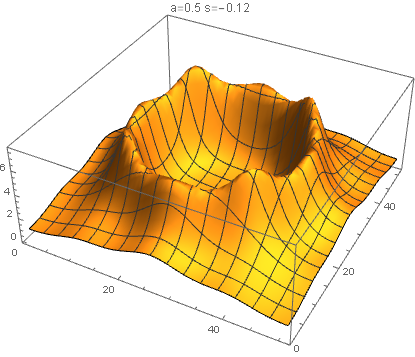}
        \caption{$t=-0.12$}
        \label{fig:ln3}
    \end{subfigure}
    \caption{Negative small imaginary time evolution of $h_{1 \bar 1}$}
%\label{fig:animals}
 %   \label{figlowneg}
\end{figure}

\begin{figure}
    \centering
    \begin{subfigure}[b]{0.2\textwidth}
        \includegraphics[width=\textwidth]{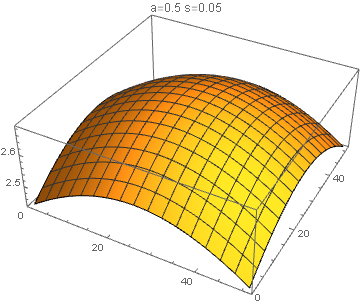}
        \caption{$t=0.05$}
        \label{fig:max1}
    \end{subfigure}
    ~ %add desired spacing between images, e. g. ~, \quad, \qquad, \hfill etc. 
    \begin{subfigure}[b]{0.2\textwidth}
        \includegraphics[width=\textwidth]{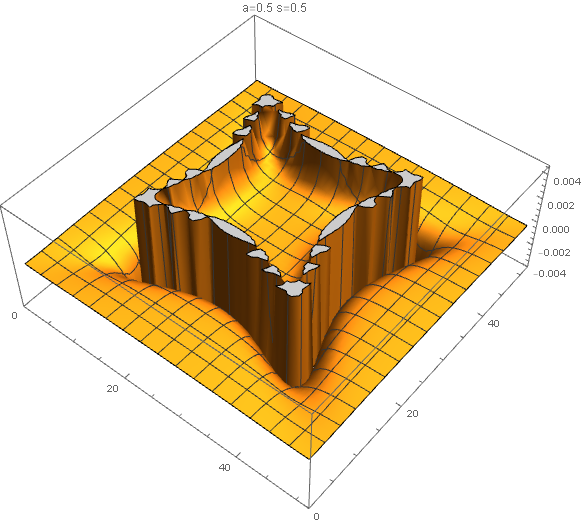}
        \caption{$t=0.5$}
        \label{fig:max2}
    \end{subfigure}
    ~ %add desired spacing between images, e. g. ~, \quad, \qquad, \hfill etc. 
    \begin{subfigure}[b]{0.2\textwidth}
        \includegraphics[width=\textwidth]{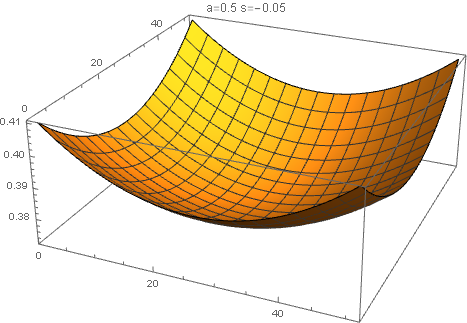}
        \caption{$t=-0.05$}
        \label{fig:max3}
    \end{subfigure}
     \begin{subfigure}[b]{0.2\textwidth}
        \includegraphics[width=\textwidth]{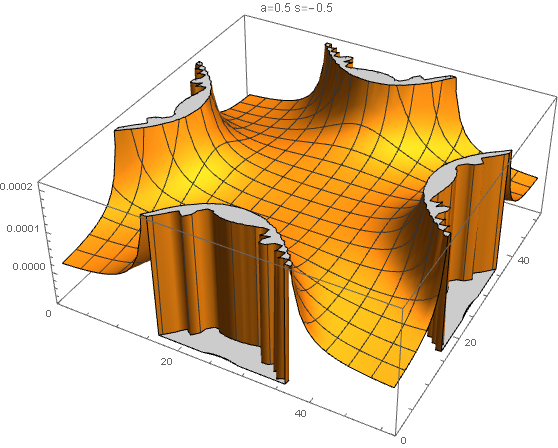}
        \caption{$t=-0.5$}
        \label{fig:max4}
    \end{subfigure}
    \caption{Evolution of $h_{1 \bar 1}$ around the maximum}
%\label{fig:animals}
    \label{figlowneg}
\end{figure}

In Figure \ref{figlowneg}, we have similar plots but for $t<0$ where the critical time is lower than but close to  
$t= 0.121$, just after Figure \ref{fig:ln3}. Due to the nature of the computation of the conformal factor, 
the transition between signs always occurs at points where the conformal factor blows up (due to zeros of the denominator in 
(\ref{e-cft2})).

\begin{figure}
    \centering
    \begin{subfigure}[b]{0.2\textwidth}
        \includegraphics[width=\textwidth]{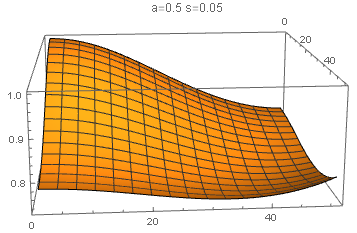}
        \caption{$t=0.3$}
        \label{fig:min1}
    \end{subfigure}
    ~ %add desired spacing between images, e. g. ~, \quad, \qquad, \hfill etc. 
    \begin{subfigure}[b]{0.2\textwidth}
        \includegraphics[width=\textwidth]{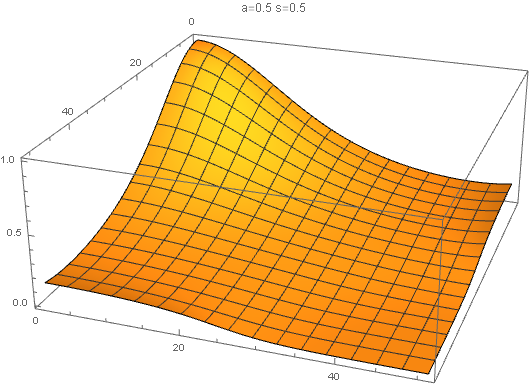}
        \caption{$t=0.5$}
        \label{fig:min2}
    \end{subfigure}
    ~ %add desired spacing between images, e. g. ~, \quad, \qquad, \hfill etc. 
    \begin{subfigure}[b]{0.2\textwidth}
        \includegraphics[width=\textwidth]{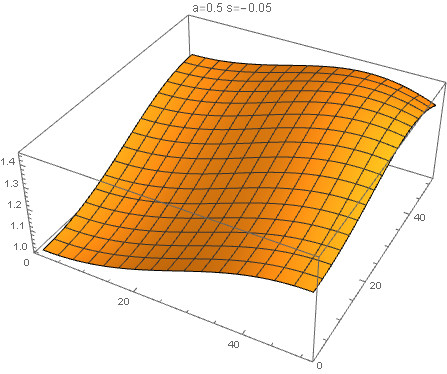}
        \caption{$t=-0.05$}
        \label{fig:min3}
    \end{subfigure}
     \begin{subfigure}[b]{0.2\textwidth}
        \includegraphics[width=\textwidth]{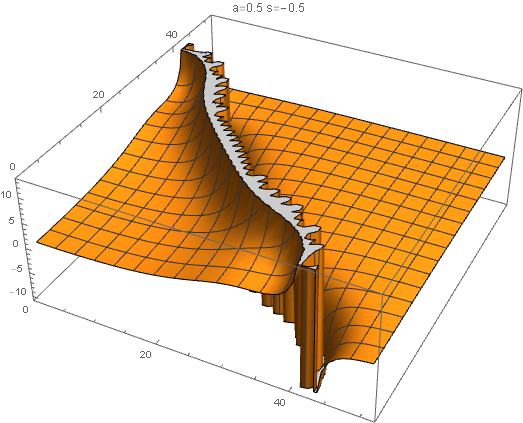}
        \caption{$t=-0.5$}
        \label{fig:min4}
    \end{subfigure}
    \caption{Evolution of $h_{1 \bar 1}$ around the minimum}
%\label{fig:animals}
 %   \label{figlowneg}
\end{figure}

Around the minimum, $H$ is close to the square of the Hamiltonian 
for the harmonic oscillator, though the region where such analogy is valid narrows with time. This phenomenon of the shrinking of the 
region of validity of the harmonic oscillator approximation is more evident for negative time evolution, where the singularity line 
approaches the origin rather quickly, as it is clear in Figure \ref{fig:min4}.

The evolution around the maximum is more complicated. Nonetheless, one also has elliptic behaviour 
around the maximum point in a (very small) region that also reduces in time. For the maximum this 
reduction is much faster and more evident since it is limited by singularity lines for both positive and negative 
imaginary time evolution.

Finally, we present the evolution of sign of the conformal factor for a value of $t$ high enough to give regions with negative conformal factor but such that we can still trust the approximation. These images, in Figure \ref{fighd}, present a very clear geometrical picture of the evolution. We also observed that taking more and more terms in the Lie series produces more and more detail 
in the patterns in the figures, but the overall structure remains similiar. 

Also, one can very clearly identify some similarities in the patterns of evolution for positive and negative time. One interesting fact is 
the similarity of the band structure of the regions of positive versus negative conformal factor around the saddle points.

%\newpage

\bigskip
\bigskip
\bigskip
\bigskip

{\bf \large{Acknowledgements:}} The authors would like to thank the referees for useful suggestions and corrections.
The authors also thank the Gulbenkian Foundation for its very successful program
``New Talents in Mathematics'' which was the driving force behind their
collaboration in this work. TR thanks also the  Gulbenkian Foundation
for his fellowship.

In addition, the authors thankfully acknowledge V. Cardoso and CENTRA/IST for computer resources, technical expertise and assistance, especially by S. Almeida. Computations were performed at the cluster ``Baltasar-Sete-S\'ois'' and supported by the H2020 ERC Consolidator Grant ``Matter and strong field gravity: New frontiers in Einstein's theory'' grant agreement no. MaGRaTh-646597

JM and JPN were partially
supported by
FCT/Portugal through the projects UID/MAT/04459/2013, PTDC/MAT-GEO/3319/2014 and by the COST Action MP1405 QSPACE.

TR thanks the Department of  Physics, Instituto Superior T\'ecnico, University of  Lisbon where most of this work was done.


\begin{thebibliography}{FMMMM}

\bibitem[BHT]{BHT} D.~Bancelin, D.~Hestroffer
and W.~Thuillot, \emph{Numerical integration of dynamical systems with Lie series}, Celest.
Mech. Dyn. Astr.  {\bf 112} (2012) 221--234.


\bibitem[BLU]{BLU} D.~Burns, E.~Lupercio and A.~Uribe, \emph{The exponential map of the complexification of {\emph Ham} in the real-analytic case},  
arXiv:1307.0493.


\bibitem[Do1]{Do1} S.K.~Donaldson, \emph{Symmetric spaces, K\"ahler geometry and Hamiltonian dynamics}, American
Math. Soc. Trans., Series 2 \textbf{196} (1999) 13--33.

\bibitem[Do2]{Do2}S.K.~Donaldson, \emph{Moment maps and diffeomorphisms}, Asian J. Math.
\textbf{3} (1999) 1--16.

\bibitem[Do3]{Do3}S.K.~Donaldson, \emph{Stability of algebraic varieties and k\"ahler geometry}, arXiv:1702.05745.

\bibitem[ED]{ED} S.~Eggl and R.~Dvorak, \emph{An introduction to common
numerical integration codes used in dynamical astronomy}, in ``Dynamics of small solar system bodies and exoplanets'', J. Souchay and R. Dvorak Eds,
Lecture Notes in Physics 790, Springer, 2010.


\bibitem[GS]{GS}   {\small 
E.-M.~Graefe, R.~Schubert, \emph{Complexified coherent states and quantum evolution with non-Hermitian Hamiltonians},
J. Phys. A  {\bf 45} (2012) 244033.}


\bibitem[Gro]{Gro} Gr\"obner, W, \emph{Die Lie-reihen und ihre Anwendungen}, Volume 3 (1967) Deutscher Verlag der Wissenschaften.


\bibitem[HLW]{HLW} E.~Hairer, C.~Lubich and G.~Wanner,
 \emph{Geometric numerical integration}, Springer, 2006.


\bibitem[Ha]{Ha}
B.C.~Hall, \emph{The {S}egal-{B}argmann ``coherent-state'' transform for Lie
  groups}, J. Funct. Anal. \textbf{122} (1994) 103--151.



\bibitem[HK]{HK}B.C.~Hall and W.D.~Kirwin, \emph{Adapted complex
structures and the geodesic flow}, Math. Ann. \textbf{350} (2011) 455--474. 


\bibitem[KMN1]{KMN1}W.D.~Kirwin, J.~Mour{\~ao}, and J.P.~Nunes,
\emph{Complex time evolution in geometric quantization and generalized
coherent state transforms}, J. Funct. Anal. Vol. \textbf{265} (2013)
1460-1493.



\bibitem[KMN2]{KMN2}W.~Kirwin, J.~Mour\~ao and J.~P. Nunes,
\emph{Complex symplectomorphisms and pseudo--K\"ahler islands in the 
quantization of toric manifolds},  Math. Ann. \textbf{364} (2016), no. 1-2, 1--28.

\bibitem[M]{M} T.~Mabuchi, \emph{Some symplectic geometry on compact K\"ahler manifolds}, 
Osaka J. Math. {\bf 24} (1987) 227--252.



\bibitem[MN]{MN}J.~Mour\~ao and J.P.~Nunes, \emph{On complexified Hamiltonian flows and 
geodesics on the space of K\"ahler metrics},  Int. Math. Res. Not. {\bf 20} (2015) 10624--10656.


\bibitem[RZ1]{RZ1}I.A~Rubinstein  and S.~Zelditch, \emph{The Cauchy problem for the homogeneous Monge-Amp\`ere equation, I. Toeplitz quantization.}, 
J. Differential Geom. {\bf 90} (2012), no. 2, 303--327. 

\bibitem[RZ2]{RZ2}I.A~Rubinstein  and S.~Zelditch, \emph{The Cauchy problem for the homogeneous Monge-Ampère equation, II. Legendre transform.}, 
Adv. Math. {\bf 228} (2011), no. 6, 2989--3025. 

\bibitem[RZ3]{RZ3}I.A~Rubinstein  and S.~Zelditch, \emph{The Cauchy problem for the homogeneous Monge-Ampère equation, III. Lifespan.}, 
J. Reine Angew. Math. {\bf 724} (2017), 105--143. 

\bibitem[Sem]{Sem}
S.~Semmes, {\it Complex Monge--Amp\`ere and symplectic manifolds},
Amer.J. Math. {\bf 114} (1992) 495--550.


\bibitem[Th]{Th} T.~Thiemann, \emph{Reality conditions inducing
transforms for quantum gauge field theory and quantum gravity},
Class.Quant.Grav. \textbf{13} (1996) 1383--1404. 

\end{thebibliography}
\end{document}